\documentclass[12pt,reqno]{amsart}
\usepackage{etex,mathtools}
\usepackage{mathdots,longtable}
\usepackage{amscd,amssymb,amsmath,multicol,tikz,float,mathrsfs}

\makeatletter
\def\subsection{\@startsection{subsection}{2}%
\z@{.5\linespacing\@plus.7\linespacing}{.5\linespacing}%
{\normalfont\bfseries}}
\makeatletter

\begin{document}

\restylefloat{table}
\newtheorem{thm}[equation]{Theorem}
\numberwithin{equation}{section}
\newtheorem{cor}[equation]{Corollary}
\newtheorem{expl}[equation]{Example}
\newtheorem{rmk}[equation]{Remark}
\newtheorem{conv}[equation]{Convention}
\newtheorem{claim}[equation]{Claim}
\newtheorem{lem}[equation]{Lemma}
\newtheorem{sublem}[equation]{Sublemma}
\newtheorem{conj}[equation]{Conjecture}
\newtheorem{defin}[equation]{Definition}
\newtheorem{diag}[equation]{Diagram}
\newtheorem{prop}[equation]{Proposition}
\newtheorem{notation}[equation]{Notation}
\newtheorem{tab}[equation]{Table}
\newtheorem{fig}[equation]{Figure}
\newcounter{bean}
\renewcommand{\theequation}{\thesection.\arabic{equation}}

\raggedbottom \voffset=-.7truein \hoffset=0truein \vsize=8truein
\hsize=6truein \textheight=8truein \textwidth=6truein
\baselineskip=18truept
\def\mapleft#1{\smash{\mathop{\longleftarrow}\limits^{#1}}}
\def\mapright#1{\ \smash{\mathop{\longrightarrow}\limits^{#1}}\ }
\def\ml#1{\,\smash{\mathop{\leftarrow}\limits^{#1}}\,}
\def\mapup#1{\Big\uparrow\rlap{$\vcenter {\hbox {$#1$}}$}}
\def\mapdown#1{\Big\downarrow\rlap{$\vcenter {\hbox {$\ssize{#1}$}}$}}
\def\mapne#1{\nearrow\rlap{$\vcenter {\hbox {$#1$}}$}}
\def\mapse#1{\searrow\rlap{$\vcenter {\hbox {$\ssize{#1}$}}$}}
\def\mapr#1{\smash{\mathop{\rightarrow}\limits^{#1}}}
\def\Mt{\widetilde{\M}}
\def\ss{\smallskip}
\def\s{\sigma}
\def\Bt{\widetilde{\mathcal{B}}}
\def\l{\lambda}
\def\Ah{\widehat{A}}
\def\Bh{\widehat{B}}
\def\vp{v_1^{-1}\pi}
\def\at{{\widetilde\alpha}}
\def\At{\widetilde{\mathcal{A}}}
\def\as{\mathscr{A}}
\def\Ast{\widetilde{\as}}
\def\Mct{\widetilde{\mathcal{M}}}
\def\sm{\wedge}
\def\la{\langle}
\def\ra{\rangle}
\def\lar{\leftarrow}
\def\ev{\text{ev}}
\def\od{\text{od}}
\def\on{\operatorname}
\def\ol#1{\overline{#1}{}}
\def\spin{\on{Spin}}
\def\cat{\on{cat}}
\def\Lbar{\overline{\Lambda}}
\def\qed{\quad\rule{8pt}{8pt}\bigskip}
\def\ssize{\scriptstyle}
\def\a{\alpha}
\def\bz{{\Bbb Z}}
\def\Rhat{\hat{R}}
\def\im{\on{im}}
\def\ct{\widetilde{C}}
\def\ext{\on{Ext}}
\def\sq{\on{Sq}}
\def\eps{\epsilon}
\def\ar#1{\stackrel {#1}{\rightarrow}}
\def\br{{\bold R}}
\def\bC{{\bold C}}
\def\bA{{\bold A}}
\def\bB{{\bold B}}
\def\bD{{\bold D}}
\def\bC{{\bold C}}
\def\bh{{\bold H}}
\def\bQ{{\bold Q}}
\def\bP{{\bold P}}
\def\bx{{\bold x}}
\def\bo{{\bold{bo}}}
\def\dh{\widehat{d}}
\def\A{\mathcal{A}}
\def\B{\mathcal{B}}
\def\si{\sigma}
\def\Vbar{{\overline V}}
\def\dbar{{\overline d}}
\def\wbar{{\overline w}}
\def\Sum{\sum}
\def\tfrac{\textstyle\frac}

\def\tb{\textstyle\binom}
\def\Si{\Sigma}
\def\w{\wedge}
\def\equ{\begin{equation}}
\def\b{\beta}
\def\G{\Gamma}
\def\L{\Lambda}
\def\g{\gamma}
\def\d{\delta}
\def\k{\kappa}
\def\psit{\widetilde{\Psi}}
\def\tht{\widetilde{\Theta}}
\def\psiu{{\underline{\Psi}}}
\def\thu{{\underline{\Theta}}}
\def\aee{A_{\text{ee}}}
\def\aeo{A_{\text{eo}}}
\def\aoo{A_{\text{oo}}}
\def\aoe{A_{\text{oe}}}
\def\vbar{{\overline v}}
\def\endeq{\end{equation}}
\def\xhat{\widehat{x}}
\def\sn{S^{2n+1}}
\def\zp{\bold Z_p}
\def\cR{{\mathcal R}}
\def\P{{\mathcal P}}
\def\cQ{{\mathcal Q}}
\def\cj{{\cal J}}
\def\zt{{\bold Z}_2}
\def\bs{{\bold s}}
\def\bof{{\bold f}}
\def\bq{{\bold Q}}
\def\be{{\bold e}}
\def\Hom{\on{Hom}}
\def\ker{\on{ker}}
\def\kot{\widetilde{KO}}
\def\coker{\on{coker}}
\def\da{\downarrow}
\def\colim{\operatornamewithlimits{colim}}
\def\zphat{\bz_2^\wedge}
\def\io{\iota}
\def\om{\omega}
\def\Prod{\prod}
\def\e{{\cal E}}
\def\zlt{\Z_{(2)}}
\def\exp{\on{exp}}
\def\abar{{\overline a}}
\def\xbar{{\overline x}}
\def\ybar{{\overline y}}
\def\zbar{{\overline z}}
\def\mbar{{\overline m}}
\def\nbar{{\overline n}}
\def\sbar{{\overline s}}
\def\kbar{{\overline k}}
\def\bbar{{\overline b}}
\def\et{{\widetilde E}}
\def\ni{\noindent}
\def\tsum{\textstyle \sum}
\def\coef{\on{coef}}
\def\den{\on{den}}
\def\lcm{\on{l.c.m.}}
\def\Ext{\operatorname{Ext}}
\def\iso{\approx}
\def\lra{\longrightarrow}
\def\vi{v_1^{-1}}
\def\ot{\otimes}
\def\psibar{{\overline\psi}}
\def\thbar{{\overline\theta}}
\def\Mh{{\widehat M}}
\def\exc{\on{exc}}
\def\ms{\medskip}
\def\ehat{{\hat e}}
\def\etao{{\eta_{\text{od}}}}
\def\etae{{\eta_{\text{ev}}}}
\def\dirlim{\operatornamewithlimits{dirlim}}
\def\gt{\widetilde{L}}
\def\lt{\widetilde{\lambda}}
\def\st{\widetilde{s}}
\def\ft{\widetilde{f}}
\def\sgd{\on{sgd}}
\def\lfl{\lfloor}
\def\rfl{\rfloor}
\def\ord{\on{ord}}
\def\gd{{\on{gd}}}
\def\rk{{{\on{rk}}_2}}
\def\nbar{{\overline{n}}}
\def\MC{\on{MC}}
\def\lg{{\on{lg}}}
\def\cH{\mathcal{H}}
\def\cS{\mathcal{S}}
\def\cP{\mathcal{P}}
\def\N{{\Bbb N}}
\def\Z{{\Bbb Z}}
\def\Q{{\Bbb Q}}
\def\R{{\Bbb R}}
\def\C{{\Bbb C}}
\def\Lb{\overline\Lambda}
\def\mo{\on{mod}}
\def\xt{\times}
\def\notimm{\not\subseteq}
\def\Remark{\noindent{\it  Remark}}
\def\kut{\widetilde{KU}}
\def\Eb{\overline E}
\def\*#1{\mathbf{#1}}
\def\0{$\*0$}
\def\1{$\*1$}
\def\22{$(\*2,\*2)$}
\def\33{$(\*3,\*3)$}
\def\ss{\smallskip}
\def\ssum{\sum\limits}
\def\dsum{\displaystyle\sum}
\def\la{\langle}
\def\ra{\rangle}
\def\on{\operatorname}
\def\proj{\on{proj}}
\def\od{\text{od}}
\def\ev{\text{ev}}
\def\o{\on{o}}
\def\U{\on{U}}
\def\lg{\on{lg}}
\def\a{\alpha}
\def\bz{{\Bbb Z}}
\def\ccM{{\Bbb M}}
\def\E{\mathcal{E}}
\def\eps{\varepsilon}
\def\bc{{\bold C}}
\def\bN{{\bold N}}
\def\bB{{\bold B}}
\def\bW{{\bold W}}
\def\nut{\widetilde{\nu}}
\def\tfrac{\textstyle\frac}
\def\b{\beta}
\def\G{\Gamma}
\def\g{\gamma}
\def\zt{{\Bbb Z}_2}
\def\zth{{\bold Z}_2^\wedge}
\def\bs{{\bold s}}
\def\bx{{\bold x}}
\def\bof{{\bold f}}
\def\bq{{\bold Q}}
\def\be{{\bold e}}
\def\lline{\rule{.6in}{.6pt}}
\def\xb{{\overline x}}
\def\xbar{{\overline x}}
\def\ybar{{\overline y}}
\def\zbar{{\overline z}}
\def\ebar{{\overline e}}
\def\nbar{{\overline n}}
\def\ubar{{\overline u}}
\def\bbar{{\overline b}}
\def\et{{\widetilde e}}
\def\M{\mathcal{M}}
\def\lf{\lfloor}
\def\rf{\rfloor}
\def\ni{\noindent}
\def\ms{\medskip}
\def\Dhat{{\widehat D}}
\def\what{{\widehat w}}
\def\Yhat{{\widehat Y}}
\def\abar{{\overline{a}}}
\def\minp{\min\nolimits'}
\def\sb{{$\ssize\bullet$}}
\def\mul{\on{mul}}
\def\N{{\Bbb N}}
\def\Z{{\Bbb Z}}
\def\bF{\mathbb{F}}
\def\Q{{\Bbb Q}}
\def\R{{\Bbb R}}
\def\C{{\Bbb C}}
\def\Xb{\overline{X}}
\def\eb{\overline{e}}
\def\notint{\cancel\cap}
\def\cS{\mathcal S}
\def\cR{\mathcal R}
\def\el{\ell}
\def\TC{\on{TC}}
\def\GC{\on{GC}}
\def\wgt{\on{wgt}}
\def\Ht{\widetilde{H}}
\def\wbar{\overline w}
\def\dstyle{\displaystyle}
\def\Sq{\on{sq}}
\def\Om{\Omega}
\def\ds{\dstyle}
\def\tz{tikzpicture}
\def\zcl{\on{zcl}}
\def\bd{\bold{d}}
\def\cM{\mathcal{M}}
\def\io{\iota}
\def\Vb#1{{\overline{V_{#1}}}}
\def\Ebar{\overline{E}}
\def\lb{\,\begin{picture}(-1,1)(1,-1)\circle*{3.5}\end{picture}\ }
\def\rlb{\,\begin{picture}(-1,1)(1,-1) \circle*{4.5}\end{picture}\ }
\def\lbb{\,\begin{picture}(-1,1)(1,-1)\circle*{8}\end{picture}\ }
\def\zp{\Z_p}
\def\lbr{\,\begin{picture}(-1,1)(1,-1)[dashed]\circle*{3.5}\end{picture}\ }
\def\llb{\,\begin{picture}(-1,1)(1,-1)\circle*{2.6}\end{picture}\ }
\def\blb{\,\begin{picture}(-1,1)(1,-1) \circle*{5.8}\end{picture}\ }
\def\alh{\widehat{\a}}
\def\bG{\bold{G}}
\setcounter{MaxMatrixCols}{15}
\title[Spin manifolds with nonzero dual Stiefel-Whitney class]
{Spin manifolds with nonzero dual Stiefel-Whitney classes of large grading}
\author{Donald M. Davis}
\address{Department of Mathematics, Lehigh University\\Bethlehem, PA 18015, USA}
\email{dmd1@lehigh.edu}
%\author{W. Stephen Wilson}
%\address{Department of Mathematics, Johns Hopkins University\\Baltimore, MD 01220, USA}
%\email{wswilsonmath@gmail.com}
\date{August 11, 2026}
\begin{abstract} The dual Stiefel-Whitney classes $\wbar_j(M)$ of a manifold $M$ are elements of $H^j(M;\zt)$ which give information about embedding and immersing $M$ in Euclidean space. We consider the question of finding, for each $n$, the largest $j$ such that there is an $n$-dimensional Spin manifold with $\wbar_j$ nonzero. We obtain upper and lower bounds. Since one theorem proves existence of manifolds without giving explicit examples, we consider a parallel question of finding explicit manifolds with $\wbar_j$ nonzero for $j$ as large as possible. \end{abstract}
\maketitle
\section{Introduction and ststement of results}\label{intro}
\subsection{$S(n)$; existence of manifolds}

Each dual Stiefel-Whitney class  $\wbar_j(M)$ of a compact differentiable\footnote{All manifolds considered in this paper will be compact and differentiable, and we will henceforth use the word ``manifold'' to mean ``compact differentiable manifold.'' All cohomology groups have coefficients in $\zt=\Z/2$.} $n$-manifold $M$ is an element of $H^j(M;\zt)$ which has the property that if $\wbar_j(M)\ne0$, then $M$ cannot be embedded in $\R^{n+j}$ nor immersed in $\R^{n+j-1}$. In \cite{DW}, the author and W.S.Wilson studied the question of finding, for each $n$, the largest $j$ such that there exists an orientable (resp.\ Spin) $n$-manifold $M$ with $\wbar_j(M)\ne0$.  This question was answered there for orientable manifolds.   Indeed, it was shown in \cite[Theorem 1.1]{DW} that the largest $j$ such that there is an orientable $n$-manifold with $\wbar_j\ne0$ is $n-\alh(n)$, where $\alh(n)$ is defined in (\ref{alh}).
\begin{equation} \label{alh}\alh(n)=\begin{cases}\a(n)&n\equiv1\pmod4\\ \a(n)+1&\text{otherwise.}\end{cases}\end{equation}
\ni Here and throughout, $\a(n)$ denotes the number of 1's in the binary expansion of $n$. However, for Spin manifolds, the question was only answered there for certain values of $n\le33$.

One aim of this paper is to extend work in \cite{DW} to find the largest $j$ such that there is an $n$-dimensional Spin manifold with nonzero $\wbar_j$. To this end, we make the following definition.
\begin{defin} We define $S(n)$ to be the largest $j$ such that there exists an $n$-dimensional Spin manifold with $\wbar_j\ne0$.
\end{defin}
 Note that $S(n)$ is a non-decreasing function of $n$, since if $M$ is an $n$-dimensional Spin manifold with $\wbar_j\ne0$, then $M\times S^1$ is $(n+1)$-dimensional with the same properties. An upper bound for $S(n)$, analogous to the sharp upper bound for orientable manifolds described above, was obtained in \cite[Theorem 1.3]{DW} using the Steenrod algebra.
 \begin{prop}$($\cite{DW}$)$\label{DWprop} We have
 $$S(n)\le\begin{cases}8k-\a(k)&\text{if }8k+1\le n\le 8k+7\\
 n-\a(n)-1&\text{if }n\equiv2^e\pmod{2^{e+2}}\\
 n-\a(n)-2&\text{if }n\equiv3\cdot 2^e\pmod{2^{e+2}}.\end{cases}$$
 \end{prop}
 \ni The similarity with the formula for orientable manifolds becomes apparent once we note that $n-\alh(n)=4k-\a(k)$ if $4k+1\le n\le4k+3$

 Using calculations of $ko$-homology groups of mod-2 Eilenberg-MacLane spaces in \cite{DW} and \cite{ko}, we can sometimes improve these upper bounds by 1.
 \begin{prop} \label{by1} If $9\le n\le12$ or $n=25$ or $n=2^e+1$ with $e\ge4$, the upper bound in Proposition \ref{DWprop} can be improved (decreased) by $1$.
 \end{prop}

 The same type of calculations can be used to prove the following lower bound for $S(n)$, which will be proved, along with Proposition \ref{by1}, in Section \ref{2.1}.
 \begin{thm} \label{thm1} For $8k\le n\le 8k+7$, we have $8k-2\a(k)\le S(n)$. This lower bound can be increased by $1$ if $13\le (n  \mod\, 16)\le 15 $
 or $18\le(n \mod\, 32)\le23$.
 \end{thm}

\subsection{$S_E(n)$; existence of Explicit manifolds}

 The proofs of the lower bound in Theorem \ref{thm1} are nonconstructive; they show that such a manifold exists, but they do not give an explicit example. The same was true for orientable manifolds in \cite[Theorem 1.1]{DW}, the result mentioned above. Because the analogous bound for all manifolds ($n-\a(n)$) was realized in \cite{Mas} by nice manifolds, namely products of $2$-power real projective spaces, the author was led in \cite{mfs} to try to find explicit orientable $n$-manifolds with $\wbar_{n-\alh(n)}\ne0$, and did so unless $n$ is a  multiple of 4 which is not a 2-power. The second aim of this paper is to try to do something like this for Spin manifolds.
 \begin{defin} Let $S_E(n)$ denote the largest $j$ such that there is a known Explicit $n$-dimensional Spin manifold with $\wbar_j\ne0$.
 \end{defin}

 Using (complex) Bott manifolds and products of 2-power quaternionic projective spaces and $\bG=G/SO(4)$, we will prove the following lower bound for $S_E(n)$. The proof of (1) is in Section \ref{2.3}, while (2) is proved in Section \ref{2.2}.
 \begin{thm}\label{thm2} \begin{enumerate}
 \item If $n\equiv2\pmod8$, then $S_E(n)\ge n-2\a(n)$.
 \item For $k\ge1$, $$S_E(8k)\ge8k-4\a(k)+\begin{cases}0&k\text{ even}\\
 2&k\text{ odd,}\end{cases}$$
 but the bound implied by (1) is larger than this if $\a(k)\ge5+\nu(k)$, where $\nu(k)$ is the exponent of $2$ in $k$.
 \end{enumerate}
 \end{thm}
 \ni This extends to give lower bounds for $S_E(n)$ for all $n\ge8$, since $S_E(n)$ is a nondecreasing function of $n$ for the same reason that $S(n)$ was.

 We will also prove the following improvement, which is optimal, in Section \ref{2.4}. This was joint work with {\tt Claude}.
 \begin{thm}\label{clthm} There is a generalized real Bott Spin manifold of dimension $13$ with $\wbar_7\ne0$. For $e\ge4$, there are generalized real Bott Spin manifolds of dimension  $2^e+4$ and $2^e+5$ with $\wbar_{2^e-1}\ne0$.\end{thm}
 \begin{cor} Let $e\ge3$ and $4\le d\le 7$. Then
 $$S(2^e+d)=S_E(2^e+d)=\begin{cases}2^e-2&e=3,\,d=4\\
 2^e-1&\text{otherwise.}\end{cases}$$\label{equal}\end{cor}

In Table \ref{tbl}, we tabulate ``$S_E$-bnd,'' our lower bound for $S_E(n)$ in Theorems \ref{thm2} and \ref{clthm}, ``$S$-bnd,'' our lower bound for $S(n)$ in Theorem \ref{thm1}, ``Upper,'' our upper bound for $S(n)$, given in Propositions \ref{DWprop} and \ref{by1}, and ``AlgUpper,'' the upper bound stated in Proposition \ref{DWprop}, derived from the Steenrod algebra. Note that the middle columns are equal for $n\le24$ and, for $e\ge4$, $n=2^e+d$ with $d\in\{0,1,4,5,6,7\}$, giving a precise value of $S(n)$ in these cases; however, the only cases when this equals $S_E(n)$ are when $8\le n\le 15$ and the values of $n$ in Corollary \ref{equal}, realized using $\bG$ and generalized real Bott manifolds. We also list manifolds with $\wbar_{S_E(n)}\ne0$. ``Bott'' refers to (complex) Bott manifolds (Theorem \ref{thm2}(1)), while ``GRBM'' refers to generalized real Bott manifolds (Theorem \ref{clthm}).

\begin{table}[H]
\centering
\caption{Bounds for $S_E(n)$ and $S(n)$}
\label{tbl}
$\begin{array}{c|cccc|c}
n&S_E\text{-bnd}&S\text{-bnd}&\text{Upper}&\text{AlgUpper}&\text{manifold}\\
\hline
8&6&6&6&6&\bG\\
9\text{-}12&6&6&6&7&\\
13\text{-}15&7&7&7&7&\text{GRBM}\\
16&12&14&14&14&HP^4\\
17&12&14&14&15&\\
18\text{-}19&14&15&15&15&\text{Bott}\\
20\text{-}23&15&15&15&15&\text{GRBM}\\
24&18&20&20&20&\bG\times HP^4\\
25&18&20&21&22&\\
26\text{-}28&20&20&22&22&\text{Bott}\\
29\text{-}31&20&21&22&22&\\
32&28&30&30&30&HP^8\\
33&28&30&30&31&\\
34\text{-}35&30&30&31&31&\text{Bott}\\
36\text{-}39&31&31&31&31&\text{GRBM}\\
40\text{-}41&34&36&38&38&\bG\times HP^8\\
42\text{-}47&36&36&38&38&\text{Bott}
\end{array}$
%\end{center}
\end{table}

Since dual Stiefel-Whitney classes are not the only way to prove nonimmersion results, the only implication that we can deduce about immersions of all Spin manifolds is given in the following result.
\begin{cor} The number $d(n)$, defined to be the smallest integer such that all $n$-dimensional Spin manifolds immerse in $\R^{d(n)}$ satisfies $d(n)\ge n +S(n)$, where for $S(n)$ we can use the lower bound in Theorem \ref{thm1} or the number in the ``$S$-\text{bnd}'' column of Table \ref{tbl}.\end{cor}

 \section{Proofs}\label{pfsec}
 \subsection{Proposition \ref{by1} and Theorem \ref{thm1}; using $ko_*$ calculations}\label{2.1}

 In \cite{DW}, we proved the following known result. Here $\chi$ is the antiautomorphism of the Steenrod algebra and $h:ko_*(X)\to H_*(X;\zt)$ is the Hurewicz homomorphism. Also $\io_k$ is the nonzero class in $H^k(K(\zt,k);\zt)$.
 \begin{prop}\label{chiprop} There exists an $n$-dimensional Spin manifold with $\wbar_{n-k}\ne0$ iff there is an element $\a\in ko_n(K(\zt,k))$ such that $\langle\chi\sq^{n-k}\io_k,h(\a)\rangle\ne0$.
 \end{prop}
\ni This uses the relationship between Spin bordism and $ko$-homology proved in \cite{ABP} to establish existence of Spin manifolds without giving explicit models.

 In \cite[Theorem 1.5]{DW}, we used detailed calculations of the Adams spectral sequence converging to $ko_*(K(\zt,k))$ for small values of $k$ to obtain all the lower bounds for $S(n)$ in the $S$-bnd column of Table \ref{tbl} for $n\le23$ and $n\in\{32,33\}$. This work also yielded the part of Proposition \ref{by1} for $9\le n\le12$.

The case $n=25$ of Proposition \ref{by1} is based on work in \cite{DW}, but was not noted there. To prove this, we must show that there does not exist an element $\a$ in $ko_{25}(K(\zt,3))$ such that $\langle\chi\sq^{22}\io_3,h_*(\a)\rangle\ne0$. Now $\chi\sq^{22}\io_3=\sq^{12}\sq^6\sq^3\sq^1\io_3$, and \cite[Figure 4,4]{DW} shows that the element of $\ext_{A_1}^{0,25}(H^*(K(\zt,3)),\zt)$ dual to $\sq^{12}\sq^6\sq^3\sq^1\io_3$ and to a certain decomposable class supports a nonzero $d_2$-differential. This eliminates the only chance of having such an element $\a$.

 In \cite{ko}, the author made a complete calculation of $ko_*(K(\zt,2))$  and used it to prove that there exists an $n$-dimensional Spin manifold with $\wbar_{n-2}\ne0$ iff $n$ is a 2-power $\ge8$. Taking the product of these manifolds corresponding to the binary expansion of $n$, using the Whitney-Cartan theorem and K\"unneth Theorem, we obtain the main part of Theorem \ref{thm1}. Taking products of these $2^e$-manifolds for $e\ge4$ with the manifolds of dimension 13 to 15 with $\wbar_7\ne0$ guaranteed by  \cite[Theorem 1.5]{DW} gives the first improvement of 1 in Theorem \ref{thm1}, and the second case of improvement follows similarly, also from \cite[Theorem 1.5]{DW}. This completes the proof of Theorem \ref{thm1}. The work in \cite{ko} also implies the $2^e+1$ part of Proposition \ref{by1}, as it says $S(2^e+1)< 2^e-1$, due to a differential in the Adams spectral sequence converging to $ko_*(K(\zt,2))$.

\subsection{Theorem \ref{thm2}(2); quaternionic projective spaces and $\bG$}\label{2.2}

Now we describe some explicit Spin manifolds.
 Our first example uses quaternionic projective space $HP^n$, which is a $4n$-dimensional Spin manifold with $H^*(HP^n)=\zt[x]/x^{n+1}$ with $|x|=4$. Its total Stiefel-Whitney class $W=\sum w_j=(1+x)^{n+1}$. Thus $\wbar_{4j}=\binom{-n-1}jx^{j}$. The formula $\binom {-a}b\equiv\binom{a+b-1}b$ mod 2 for $a,b\ge0$ is useful. If $n=2^e$, then $\wbar_{4n-4}(HP^n)=\binom{-2^e-1}{2^e-1}x^{n-1}$, with coefficient equal to $\binom{2^{e+1}-1}{2^e-1}\equiv1$ by Lucas's Theorem. Thus $HP^n$ has $\wbar_{4n-4}\ne0$ if $n$ is a 2-power. If $n=\sum 2^{e_i}$ for distinct $e_i\ge3$,\footnote{We could include $e_i=2$, but $HP^1$ would perform the same function as the $(S^1)^4$ that we have been using.} then $\prod HP^{2^{e_i-2}}$ has $\wbar_{n-4\a(n)}\ne0$. These, producted with $(S^1)^t$ for $1\le t\le 7$, give the case of Theorem \ref{thm2}(2) with $k$ even.

  It was pointed out to the author by {\tt ChatGPT} that the symmetric space $\bG=G_2/SO(4)$ is an 8-dimensional Spin manifold with $\wbar_6\ne0$. This was shown long ago in \cite{BH}. This is optimal; we use $\bG$ as a replacement for $HP^2$. By using products of distinct $HP^{2^e}$'s for $e\ge2$, $(S^1)^t$ for $1\le t\le7$, and $\bG$, we obtain Theorem \ref{thm2}(2).

\subsection{Theorem \ref{thm2}(1); (complex) Bott manifolds}\label{2.3}

  In \cite{mfs}, the author used real Bott manifolds to obtain orientable $n$-manifolds with $\wbar_{n-\a(n)}\ne0$ for $n\equiv1$ mod 4. The complex analogue is a $2n$-dimensional Spin manifold with $\wbar_{2n-2\a(n)}\ne0$. These are our examples implying Theorem \ref{thm2}(1). We elaborate.

  Let $\bF=\R$ or $\C$, and $d=1$ (resp.\ 2) if $\bF=\R$ (resp.\ $\C$).
  An $\bF$-Bott tower is a sequence of fibrations
  $$B_n\to\cdots \to B_1\to B_0=*$$
  such that $B_{j+1}=P(\eps_{\bF}\oplus\zeta_j)$ with $\zeta_j$ an $\bF$-line bundle over $B_j$, $\eps_{\bF}$ a trivial $\bF$-line bundle, and $P$  the projectification. Then $B_n$ is a $dn$-manifold, and there are canonical classes $x_i$ in $H^d(B_j)$ for $i\le j$ and
  $$w_d(\zeta_j)=\sum_{i< j}a_{i,j}x_j$$
  for some $a_{i,j}\in\zt$. These $a_{i,j}$ form an $n$-by-$n$ strictly upper-triangular binary matrix $A$ which determines $B_n$. Then, from (\cite{CS}, \cite{Ds}), we have
  \begin{equation}\label{coh}H^*(B_n)=\zt[x_1,\ldots,x_n]/(x_j^2=\sum_{i< j}a_{i,j}x_ix_j)\end{equation}
  and
  $$w_{dk}(\tau(B_n))=\si_k(\{\sum_{i<j}a_{i,j}x_i,\ 2\le j\le n\}),$$
  with $\si_k$ the $k$th elementary symmetric function, and $\tau$ the tangent bundle.

  The author proved in \cite{mfs} that, for $\bF=\R$, for the matrix $A$ with
  $$a_{i,j}=\begin{cases}1&j=i+1\text{ or }n,\ i\le n-2\\
  0&\text{otherwise,}\end{cases}$$
  the corresponding real Bott matrix is orientable and has $\wbar_{n-\a(n)}\ne0$ if $n\equiv1$ mod 4. %Since $\a(4k+1)=\alh(4k+1)=\alh(4k+2)=\alh(4k+3)$, producting with one or two copies of $S^1$ yields manifolds with the same property.
  If we form the complex Bott manifold of (real) dimension $2n$ with the same matrix $A$, the same proof yields $w_2=0$ and $\wbar_{2n-2\a(n)}\ne0$ if $n\equiv1$ mod 4. This proves Theorem \ref{thm2}(1), with the $n$ in the theorem corresponding to the $2n$ in the preceding sentence.

\subsection{Theorem \ref{clthm}; generalized real Bott manifolds}\label{2.4}

We close this section by proving Theorem \ref{clthm}. Rather than recapitulating the general case of a generalized real Bott manifold, which involves lots of subscripts and is described in \cite{DU} and \cite{Turk}, we describe directly the case $n=2^e+4$, $e\ge4$ in our theorem. Let $(n_1,n_2,n_3,n_4)=(2,2,2^e-3,3)$. There is a tower of fibrations
$$B_4\to B_3\to B_2\to B_1\to B_0=*$$
in which $B_{j+1}=P(\eps\oplus V_j)$, where $V_j$ is a sum of $n_j$ line bundles over $B_j$. (The space $B_4$ can be considered as a small cover over $\Delta^{n_1}\times\Delta^{n_2}\times\Delta^{n_3}\times\Delta^{n_4}$, but that is not our approach.)

The partitioned matrix $M=[M_1|M_2|M_3|M_4]=$

\begin{table}[H]
$\left[\begin{array}{cc|cc|ccccccccccc|ccc}
\centering
1&1&0&0&1&1&0&0&0&0&0&1&0&\cdots&0&1&0&1\\
0&0&1&1&1&1&1&1&1&1&1&0&0&\cdots&0&1&0&1\\
0&0&0&0&1&1&1&1&1&1&1&1&1&\cdots&1&0&1&1\\
0&0&0&0&0&0&0&0&0&0&0&0&0&\cdots&0&1&1&1
\end{array}\right],$\end{table}

\ni discovered by {\tt Claude},
describes the tower completely. The four ellipses ($\cdots$) refer to sequences of $(2^e-13)$ 0's, 0's, 1's, and 0's, resp. There are canonical elements $x_i\in H^1(B_j)$ for $1\le i\le j\le4$. Let $\bx=(x_1,x_2,x_3,x_4)^T$. The dot product of  the $t$th column of $M_k$ with $\bx$ equals $w_1(L_{k,t})+x_k$, where $L_{k,t}$ is the $t$th line bundle in $V_k$. For example, $w_1(L_{4,1})=x_1+x_2$. Then
$$H^*(B_4)=\zt[x_1,x_2,x_3,x_4]/(R_1,R_2,R_3,R_4),$$
where $R_k$ is the relation
$$x_k^{n_k+1}=\sum_{i=1}^{n_k}\si_i(S_k),$$
with $S_k$ the set of the $n_k$ expressions $w_1(L_{k,t})$, and $\si_i$ is the $i$th elementary symmetric polynomial.  Explicitly, we have

\begin{align}\nonumber x_1^3=&\,0\\
x_2^3=&\,0\nonumber\\
x_3^{2^e-2}=&\,(x_1+x_2)x_3^{2^e-3}+(x_1^2+x_2^2+x_1x_2)x_3^{2^e-4}+(x_1^3+x_1^2x_2+x_1x_2^2+x_2^3)x_3^{2^e-5}\nonumber\\
&\,+(x_1^3x_2+x_1x_2^3+x_2^4)x_3^{2^e-6}+(x_1x_2^4+x_2^5)x_3^{2^e-7}+(x_1^2x_2^4+x_1x_2^5+x_2^6)x_3^{2^e-8}\nonumber\\
&\,+(x_1^3x_2^4+x_1^2x_2^5+x_1x_2^6+x_2^7)x_3^{2^e-9}+(x_1^3x_2^5+x_1x_2^7)x_3^{2^e-10}\label{3}\\
x_4^4=&\,(x_1^2+x_2^2+x_3^2+x_1x_3+x_2x_3)x_4^2+(x_1^2x_3+x_2^2x_3+x_1x_3^2+x_2x_3^2)x_4.\label{4}
\end{align}
\ni For example, the coefficient of $x_4^2$ in $x_4^4$ is $\si_2(x_1+x_2,x_3,x_1+x_2+x_3)$. A basis for $H^*(B_4)$ is $\{x_1^{i_1}x_2^{i_2}x_3^{i_3}x_4^{i_4}:\,0\le i_t\le n_k\}$.

We use the standard result (\cite[\S3]{MP}) that $\wbar_k$ is nonzero in an $n$-manifold  $M$ if and only if there is an $x\in H^{n-k}(M)$ with $\chi\sq^k(x)\ne0$.
 Using $(a,b,c,d)$ to represent $x_1^ax_2^bx_3^cx_4^d$, one easily shows
\begin{gather}\chi\sq^{2^e-1}\label{chi}(x_2x_3x_4^3)=\sq^{2^{e-1}}\sq^{2^{e-2}}\cdots\sq^2\sq^1(x_2x_3x_4^3)\\=(0,2,1,2^e+1)+(0,1,2,2^e+1)\nonumber+(0,2,2,2^e)+(0,1,2^e,3)+(0,2^e,1,3)\end{gather}
for $e\ge2$.
Showing that this is nonzero when the relations are applied was first done for $e=4$ and $5$ by {\tt Claude} in the discovery process, and then verified by the author for $e=4$ using a {\tt Maple} program. We can prove it for all $e$ as follows.

First we use the relation (\ref{4}) for $x_4^4$ together with $x_1^3=x_2^3=0$ to reduce the exponent of $x_4$ to values $\le3$. We claim that the terms reduce as follows.
\begin{align}(0,2,1,2^e+1)&=(2,2,2^e-2,2)+(2,2,2^e-3,3)
+(1,2,2^e-1,2)\nonumber\\ &\quad+(1,2,2^e-2,3)
+(0,2,2^e-1,3)\nonumber\\
(0,1,2,2^e+1)&=(2,2,2^e-2,2)+(2,2,2^e-3,3)+(2,1,2^e-1,2)\nonumber\\&\quad+(2,1,2^e-2,3)+(1,1,2^e,2)
 +(1,1,2^e-1,3)\nonumber\\&\quad+(0,2,2^e,2)
 +(0,2,2^e-1,3)+(0,1,2^e,3)\nonumber\\
(0,2,2,2^e)&=(2,2,2^e-1,1)+(2,2,2^e-2,2)+(1,2,2^e,1)\nonumber\\
&\quad+(1,2,2^e-1,2)+(0,2,2^e,2)\label{5}\\
(0,1,2^e,3)&=(0,1,2^e,3)\nonumber\\
(0,2^e,1,3)&=0.\nonumber
\end{align}

We now prove this claim. Write the relation (\ref{4}) as $x_4^4=Ax_4^2+Bx_4$. To find a formula for the coefficients making $x_4^n=c_{1,n}x_4^3+c_{2,n}x_4^2+c_{3,n}x_4$, multiply by $x_4$ and reduce using $x_4^4=Ax_4^2+Bx_4$. The initial condition for this recursion is $(c_{1,4},c_{2,4},c_{3,4})=(0,A,B)$. The recursion is
$$(c_{1,n+1},c_{2,n+1},c_{3,n+1})=(c_{2,n},c_{3,n}+Ac_{1,n},Bc_{1,n})
$$ and the solution is
\begin{align*} &(c_{1,n},c_{2,n},c_{3,n})\\
=&\biggl(\sum_{2i+3j=n-3}\binom{i+j}iA^iB^j,\sum_{2i+3j=n-2}\binom{i+j}iA^iB^j,\sum_{2i+3j=n-1}\binom{i+j-1}iA^iB^j\biggr).\end{align*}

We demonstrate the reduction  of $x_2^2x_3^2x_4^{2^e}$ which results in (\ref{5}); the two equations preceding (\ref{5}) are obtained similarly. One shows easily that $c_{1,2^e}=0$, $c_{2,2^e}\equiv A^{2^{e-1}-1}+A^{2^{e-1}-4}B^2
\mod \ B^3$, and $c_{3,2^e}\equiv A^{2^{e-1}-2}B\mod\ B^2$. Since we are reducing $x_2^2x_3^2x_4^{2^e}$, terms divisible by $x_2$ or by $x_1^3$ can be omitted. We obtain
\begin{align*}&x_2^2x_3^2x_4^{2^e}\\
=&x_2^2x_3^2\biggl(x_4^2\bigl((x_1^2+x_3^2+x_1x_3)^{2^{e-1}-1}+(x_3^2)^{2^{e-1}-4}(x_1x_3^2)^2\bigr)\\
&\qquad +x_4(x_3^2)^{2^{e-1}-2}(x_1^2x_3+x_1x_3^2)\biggr)\\
=&x_2^2x_3^2\bigl(x_4^2(x_3^{2^e-2}+x_1x_3^{2^e-3}+x_1^2x_3^{2^e-4})+x_4x_3^{2^e-4}(x_1^2x_3+x_1x_3^2)\bigr),\end{align*}
which are the five terms listed in (\ref{5}).

Next we use (\ref{3}) and $x_1^3=x_2^3=0$ to see which of the monomials on the right hand side of the equations surrounding (\ref{5}) equal the top class $x_1^2x_2^2x_3^{2^e-3}x_4^3$, and which equal 0. We claim that all except
$(2,2,2^e-3,3)$, $(1,2,2^e-2,3)$, $
(2,1,2^e-2,3)$, and $(1,1,2^e-1,3)$ reduce to 0, while these reduce to $(2,2,2^e-3,3)$. This implies that only the second of the five equations and hence the second of the five terms in (\ref{chi}) is nonzero.
We conclude that $\wbar_{2^e-1}\ne0$ since $\chi\sq^{2^e-1}(x_2x_3x_4^3)\ne0$.

We now prove this claim. Since the $x_4^k$ factors are not affected by this reduction, only the monomials ending in $x_4^3$ have a chance to be nonzero. We illustrate the reduction using $x_1x_2x_3^{2^e-1}x_4^3$. Multiply the relation (\ref{3}) by $x_3$ and ignore all terms divisible by $x_1^2$ or $x_2^2$. We obtain

\begin{align*}x_1x_2x_3^{2^e-1}x_4^3&=x_1x_2\bigl((x_1+x_2)x_3^{2^e-2}+x_1x_2x_3^{2^e-3}\bigr)x_4^3\\
&=\bigl((x_1^2x_2+x_1x_2^2)(x_1+x_2)x_3^{2^e-3}+x_1^2x_2^2x_3^{2^e-3}\bigr)x_4^2\\
&=x_1^2x_2^2x_3^{2^e-3}x_4^3\ne0.\end{align*}
Other monomials are handled similarly.

We can obtain $n=13$
part of Theorem \ref{clthm} similarly using the matrix

\begin{table}[H]
$\left[\begin{array}{cc|ccc|ccccc|ccc}
\centering
1&1&1&1&1&1&1&1&1&0&1&0&1\\
0&0&1&1&1&1&0&0&0&1&1&0&1\\
0&0&0&0&0&1&1&1&1&1&1&0&1\\
0&0&0&0&0&0&0&0&0&0&1&1&1
\end{array}\right]$\end{table}

\ni and $\chi\sq^7(x_1x_2x_3x_4^3)$. Again this was discovered by {\tt Claude}. We can generalize this to all $2^e+5$ by appending 0's and 1's to the third block as we did for the $2^e+4$ case. This is not really necessary since, for $e\ge4$, examples of $(2^e+5)$-dimensional Spin manifolds with $\wbar_{2^e-1}\ne0$ can be obtained by crossing the $(2^4+4)$-dimensional example with $S^1$.

\section{Other manifolds}
  We tested the generalized Dold manifolds of \cite{SZ} and found that there are some values of $n$
  for which their total dual Stiefel-Whitney class has top nonzero grading equal to our lower bounds for $S_E(n)$ in Theorem \ref{thm2}, but they never improve upon it, so we will not provide details.

We studied other  generalized real Bott manifolds in addition to those in the dimensions mentioned above. Based on many trials by {\tt Claude}
we feel that it is unlikely that there exist generalized real Bott Spin manifolds of dimension 16 or 17 with $\wbar_{12}$ or $\wbar_{13}$ nonzero, or of dimension 18 or 19 with $\wbar_{15}$ nonzero, or of dimension 24 or 25 with $\wbar_{18}$ or $\wbar_{19}$ nonzero.

Let $B_n$ be a real Bott $n$-manifold with matrix $A$ having rows $R_i$. It
 is shown in \cite{Ds} that $B_n$ is Spin if $R_i\cdot R_i\equiv0$ mod 2 for all $i$ and, if $\widehat{R}_i$ is obtained from $R_i$ by letting $a_{i,i}$ equal the mod 2 reduction of $\frac12 R_i\cdot R_i$,
 then $R_j\cdot \widehat{R}_i\equiv0$ mod 2 whenever $j<i$. {\tt Claude} found a 14-dimensional real Bott Spin manifold with $\wbar_7\ne0$, but we feel that this was superceded by the generalized real Bott Spin manifolds described earlier in the paper.

\def\line{\rule{.6in}{.6pt}}

\end{document}